\overfullrule=0pt
\magnification=1200
\hsize=12.5cm    
\vsize=19cm
\hoffset=1cm
\font\grand=cmr10 at 14pt

\overfullrule=0pt
\magnification=1200
\hsize=11.25cm    
\vsize=19cm
\hoffset=1cm

\font\grand=cmr10 at 12pt

\def\N{\noindent}
\def\S{\smallskip \par}
\def\M{\medskip \par}
\def\B{\bigskip \par}
\def\BB{\bigskip \bigskip \par}
\def\P{\noindent {\it Proof.} }

\catcode`@=11
\def\@height{height}
\def\@depth{depth}
\def\@width{width}

\newcount\@tempcnta
\newcount\@tempcntb

\newdimen\@tempdima
\newdimen\@tempdimb

\newbox\@tempboxa

\def\@ifnextchar#1#2#3{\let\@tempe #1\def\@tempa{#2}\def\@tempb{#3}\futurelet
    \@tempc\@ifnch}
\def\@ifnch{\ifx \@tempc \@sptoken \let\@tempd\@xifnch
      \else \ifx \@tempc \@tempe\let\@tempd\@tempa\else\let\@tempd\@tempb\fi
      \fi \@tempd}
\def\@ifstar#1#2{\@ifnextchar *{\def\@tempa*{#1}\@tempa}{#2}}

\def\@whilenoop#1{}
\def\@whilenum#1\do #2{\ifnum #1\relax #2\relax\@iwhilenum{#1\relax 
     #2\relax}\fi}
\def\@iwhilenum#1{\ifnum #1\let\@nextwhile=\@iwhilenum 
         \else\let\@nextwhile=\@whilenoop\fi\@nextwhile{#1}}

\def\@whiledim#1\do #2{\ifdim #1\relax#2\@iwhiledim{#1\relax#2}\fi}
\def\@iwhiledim#1{\ifdim #1\let\@nextwhile=\@iwhiledim 
        \else\let\@nextwhile=\@whilenoop\fi\@nextwhile{#1}}

\newdimen\@wholewidth
\newdimen\@halfwidth
\newdimen\unitlength \unitlength =1pt
\newbox\@picbox
\newdimen\@picht

\def\@nnil{\@nil}
\def\@empty{}
\def\@fornoop#1\@@#2#3{}

\def\@for#1:=#2\do#3{\edef\@fortmp{#2}\ifx\@fortmp\@empty \else
    \expandafter\@forloop#2,\@nil,\@nil\@@#1{#3}\fi}

\def\@forloop#1,#2,#3\@@#4#5{\def#4{#1}\ifx #4\@nnil \else
       #5\def#4{#2}\ifx #4\@nnil \else#5\@iforloop #3\@@#4{#5}\fi\fi}

\def\@iforloop#1,#2\@@#3#4{\def#3{#1}\ifx #3\@nnil 
       \let\@nextwhile=\@fornoop \else
      #4\relax\let\@nextwhile=\@iforloop\fi\@nextwhile#2\@@#3{#4}}

\def\@tfor#1:=#2\do#3{\xdef\@fortmp{#2}\ifx\@fortmp\@empty \else
    \@tforloop#2\@nil\@nil\@@#1{#3}\fi}
\def\@tforloop#1#2\@@#3#4{\def#3{#1}\ifx #3\@nnil 
       \let\@nextwhile=\@fornoop \else
      #4\relax\let\@nextwhile=\@tforloop\fi\@nextwhile#2\@@#3{#4}}

\def\@makepicbox(#1,#2){\leavevmode\@ifnextchar 
   [{\@imakepicbox(#1,#2)}{\@imakepicbox(#1,#2)[]}}

\long\def\@imakepicbox(#1,#2)[#3]#4{\vbox to#2\unitlength
   {\let\mb@b\vss \let\mb@l\hss\let\mb@r\hss
    \let\mb@t\vss
    \@tfor\@tempa :=#3\do{\expandafter\let
        \csname mb@\@tempa\endcsname\relax}%
\mb@t\hbox to #1\unitlength{\mb@l #4\mb@r}\mb@b}}

\def\picture(#1,#2){\@ifnextchar({\@picture(#1,#2)}{\@picture(#1,#2)(0,0)}}

\def\@picture(#1,#2)(#3,#4){\@picht #2\unitlength
\setbox\@picbox\hbox to #1\unitlength\bgroup 
\hskip -#3\unitlength \lower #4\unitlength \hbox\bgroup\ignorespaces}

\def\endpicture{\egroup\hss\egroup\ht\@picbox\@picht
\dp\@picbox\z@\leavevmode\box\@picbox}

\long\def\put(#1,#2)#3{\@killglue\raise#2\unitlength\hbox to \z@{\kern
#1\unitlength #3\hss}\ignorespaces}

\long\def\multiput(#1,#2)(#3,#4)#5#6{\@killglue\@multicnt=#5\relax
\@xdim=#1\unitlength
\@ydim=#2\unitlength
\@whilenum \@multicnt > 0\do
{\raise\@ydim\hbox to \z@{\kern
\@xdim #6\hss}\advance\@multicnt \m@ne\advance\@xdim
#3\unitlength\advance\@ydim #4\unitlength}\ignorespaces}

\def\@killglue{\unskip\@whiledim \lastskip >\z@\do{\unskip}}

\def\thinlines{\let\@linefnt\tenln \let\@circlefnt\tencirc
  \@wholewidth\fontdimen8\tenln \@halfwidth .5\@wholewidth}
\def\thicklines{\let\@linefnt\tenlnw \let\@circlefnt\tencircw
  \@wholewidth\fontdimen8\tenlnw \@halfwidth .5\@wholewidth}

\def\linethickness#1{\@wholewidth #1\relax \@halfwidth .5\@wholewidth}

\def\shortstack{\@ifnextchar[{\@shortstack}{\@shortstack[c]}}

\def\@shortstack[#1]{\leavevmode
\vbox\bgroup\baselineskip-1pt\lineskip 3pt\let\mb@l\hss
\let\mb@r\hss \expandafter\let\csname mb@#1\endcsname\relax
\let\\\@stackcr\@ishortstack}

\def\@ishortstack#1{\halign{\mb@l ##\unskip\mb@r\cr #1\crcr}\egroup}

\def\@stackcr{\@ifstar{\@ixstackcr}{\@ixstackcr}}
\def\@ixstackcr{\@ifnextchar[{\@istackcr}{\cr\ignorespaces}}

\def\@istackcr[#1]{\cr\noalign{\vskip #1}\ignorespaces}

\newif\if@negarg

\def\droite(#1,#2)#3{\@xarg #1\relax \@yarg #2\relax
\@linelen=#3\unitlength
\ifnum\@xarg =0 \@vline 
  \else \ifnum\@yarg =0 \@hline \else \@sline\fi
\fi}

\def\@sline{\ifnum\@xarg< 0 \@negargtrue \@xarg -\@xarg \@yyarg -\@yarg
  \else \@negargfalse \@yyarg \@yarg \fi
\ifnum \@yyarg >0 \@tempcnta\@yyarg \else \@tempcnta -\@yyarg \fi
\ifnum\@tempcnta>6 \@badlinearg\@tempcnta0 \fi
\ifnum\@xarg>6 \@badlinearg\@xarg 1 \fi
\setbox\@linechar\hbox{\@linefnt\@getlinechar(\@xarg,\@yyarg)}%
\ifnum \@yarg >0 \let\@upordown\raise \@clnht\z@
   \else\let\@upordown\lower \@clnht \ht\@linechar\fi
\@clnwd=\wd\@linechar
\if@negarg \hskip -\wd\@linechar \def\@tempa{\hskip -2\wd\@linechar}\else
     \let\@tempa\relax \fi
\@whiledim \@clnwd <\@linelen \do
  {\@upordown\@clnht\copy\@linechar
   \@tempa
   \advance\@clnht \ht\@linechar
   \advance\@clnwd \wd\@linechar}%
\advance\@clnht -\ht\@linechar
\advance\@clnwd -\wd\@linechar
\@tempdima\@linelen\advance\@tempdima -\@clnwd
\@tempdimb\@tempdima\advance\@tempdimb -\wd\@linechar
\if@negarg \hskip -\@tempdimb \else \hskip \@tempdimb \fi
\multiply\@tempdima \@m
\@tempcnta \@tempdima \@tempdima \wd\@linechar \divide\@tempcnta \@tempdima
\@tempdima \ht\@linechar \multiply\@tempdima \@tempcnta
\divide\@tempdima \@m
\advance\@clnht \@tempdima
\ifdim \@linelen <\wd\@linechar
   \hskip \wd\@linechar
  \else\@upordown\@clnht\copy\@linechar\fi}

\def\@hline{\ifnum \@xarg <0 \hskip -\@linelen \fi
\vrule \@height \@halfwidth \@depth \@halfwidth \@width \@linelen
\ifnum \@xarg <0 \hskip -\@linelen \fi}

\def\@getlinechar(#1,#2){\@tempcnta#1\relax\multiply\@tempcnta 8
\advance\@tempcnta -9 \ifnum #2>0 \advance\@tempcnta #2\relax\else
\advance\@tempcnta -#2\relax\advance\@tempcnta 64 \fi
\char\@tempcnta}

\def\vector(#1,#2)#3{\@xarg #1\relax \@yarg #2\relax
\@tempcnta \ifnum\@xarg<0 -\@xarg\else\@xarg\fi
\ifnum\@tempcnta<5\relax
\@linelen=#3\unitlength
\ifnum\@xarg =0 \@vvector 
  \else \ifnum\@yarg =0 \@hvector \else \@svector\fi
\fi
\else\@badlinearg\fi}

\def\@hvector{\@hline\hbox to 0pt{\@linefnt 
\ifnum \@xarg <0 \@getlarrow(1,0)\hss\else
    \hss\@getrarrow(1,0)\fi}}

\def\@vvector{\ifnum \@yarg <0 \@downvector \else \@upvector \fi}

\def\@svector{\@sline
\@tempcnta\@yarg \ifnum\@tempcnta <0 \@tempcnta=-\@tempcnta\fi
\ifnum\@tempcnta <5
  \hskip -\wd\@linechar
  \@upordown\@clnht \hbox{\@linefnt  \if@negarg 
  \@getlarrow(\@xarg,\@yyarg) \else \@getrarrow(\@xarg,\@yyarg) \fi}%
\else\@badlinearg\fi}

\def\@getlarrow(#1,#2){\ifnum #2 =\z@ \@tempcnta='33\else
\@tempcnta=#1\relax\multiply\@tempcnta \sixt@@n \advance\@tempcnta
-9 \@tempcntb=#2\relax\multiply\@tempcntb \tw@
\ifnum \@tempcntb >0 \advance\@tempcnta \@tempcntb\relax
\else\advance\@tempcnta -\@tempcntb\advance\@tempcnta 64
\fi\fi\char\@tempcnta}

\def\@getrarrow(#1,#2){\@tempcntb=#2\relax
\ifnum\@tempcntb < 0 \@tempcntb=-\@tempcntb\relax\fi
\ifcase \@tempcntb\relax \@tempcnta='55 \or 
\ifnum #1<3 \@tempcnta=#1\relax\multiply\@tempcnta
24 \advance\@tempcnta -6 \else \ifnum #1=3 \@tempcnta=49
\else\@tempcnta=58 \fi\fi\or 
\ifnum #1<3 \@tempcnta=#1\relax\multiply\@tempcnta
24 \advance\@tempcnta -3 \else \@tempcnta=51\fi\or 
\@tempcnta=#1\relax\multiply\@tempcnta
\sixt@@n \advance\@tempcnta -\tw@ \else
\@tempcnta=#1\relax\multiply\@tempcnta
\sixt@@n \advance\@tempcnta 7 \fi\ifnum #2<0 \advance\@tempcnta 64 \fi
\char\@tempcnta}

\def\@vline{\ifnum \@yarg <0 \@downline \else \@upline\fi}

\def\@upline{\hbox to \z@{\hskip -\@halfwidth \vrule \@width \@wholewidth
   \@height \@linelen \@depth \z@\hss}}

\def\@downline{\hbox to \z@{\hskip -\@halfwidth \vrule \@width \@wholewidth
   \@height \z@ \@depth \@linelen \hss}}

\def\@upvector{\@upline\setbox\@tempboxa\hbox{\@linefnt\char'66}\raise 
     \@linelen \hbox to\z@{\lower \ht\@tempboxa\box\@tempboxa\hss}}

\def\@downvector{\@downline\lower \@linelen
      \hbox to \z@{\@linefnt\char'77\hss}}

\def\dashbox#1(#2,#3){\leavevmode\hbox to \z@{\baselineskip \z@%
\lineskip \z@%
\@dashdim=#2\unitlength%
\@dashcnt=\@dashdim \advance\@dashcnt 200
\@dashdim=#1\unitlength\divide\@dashcnt \@dashdim
\ifodd\@dashcnt\@dashdim=\z@%
\advance\@dashcnt \@ne \divide\@dashcnt \tw@ 
\else \divide\@dashdim \tw@ \divide\@dashcnt \tw@
\advance\@dashcnt \m@ne
\setbox\@dashbox=\hbox{\vrule \@height \@halfwidth \@depth \@halfwidth
\@width \@dashdim}\put(0,0){\copy\@dashbox}%
\put(0,#3){\copy\@dashbox}%
\put(#2,0){\hskip-\@dashdim\copy\@dashbox}%
\put(#2,#3){\hskip-\@dashdim\box\@dashbox}%
\multiply\@dashdim 3 
\fi
\setbox\@dashbox=\hbox{\vrule \@height \@halfwidth \@depth \@halfwidth
\@width #1\unitlength\hskip #1\unitlength}\@tempcnta=0
\put(0,0){\hskip\@dashdim \@whilenum \@tempcnta <\@dashcnt
\do{\copy\@dashbox\advance\@tempcnta \@ne }}\@tempcnta=0
\put(0,#3){\hskip\@dashdim \@whilenum \@tempcnta <\@dashcnt
\do{\copy\@dashbox\advance\@tempcnta \@ne }}%
\@dashdim=#3\unitlength%
\@dashcnt=\@dashdim \advance\@dashcnt 200
\@dashdim=#1\unitlength\divide\@dashcnt \@dashdim
\ifodd\@dashcnt \@dashdim=\z@%
\advance\@dashcnt \@ne \divide\@dashcnt \tw@
\else
\divide\@dashdim \tw@ \divide\@dashcnt \tw@
\advance\@dashcnt \m@ne
\setbox\@dashbox\hbox{\hskip -\@halfwidth
\vrule \@width \@wholewidth 
\@height \@dashdim}\put(0,0){\copy\@dashbox}%
\put(#2,0){\copy\@dashbox}%
\put(0,#3){\lower\@dashdim\copy\@dashbox}%
\put(#2,#3){\lower\@dashdim\copy\@dashbox}%
\multiply\@dashdim 3
\fi
\setbox\@dashbox\hbox{\vrule \@width \@wholewidth 
\@height #1\unitlength}\@tempcnta0
\put(0,0){\hskip -\@halfwidth \vbox{\@whilenum \@tempcnta < \@dashcnt
\do{\vskip #1\unitlength\copy\@dashbox\advance\@tempcnta \@ne }%
\vskip\@dashdim}}\@tempcnta0
\put(#2,0){\hskip -\@halfwidth \vbox{\@whilenum \@tempcnta< \@dashcnt
\relax\do{\vskip #1\unitlength\copy\@dashbox\advance\@tempcnta \@ne }%
\vskip\@dashdim}}}\@makepicbox(#2,#3)}

\newif\if@ovt 
\newif\if@ovb 
\newif\if@ovl 
\newif\if@ovr 
\newdimen\@ovxx
\newdimen\@ovyy
\newdimen\@ovdx
\newdimen\@ovdy
\newdimen\@ovro
\newdimen\@ovri

\def\@getcirc#1{\@tempdima #1\relax \advance\@tempdima 2pt\relax
  \@tempcnta\@tempdima
  \@tempdima 4pt\relax \divide\@tempcnta\@tempdima
  \ifnum \@tempcnta > 10\relax \@tempcnta 10\relax\fi
  \ifnum \@tempcnta >\z@ \advance\@tempcnta\m@ne
    \else \@warning{Oval too small}\fi
  \multiply\@tempcnta 4\relax
  \setbox \@tempboxa \hbox{\@circlefnt
  \char \@tempcnta}\@tempdima \wd \@tempboxa}

\def\@put#1#2#3{\raise #2\hbox to \z@{\hskip #1#3\hss}}

\def\oval(#1,#2){\@ifnextchar[{\@oval(#1,#2)}{\@oval(#1,#2)[]}}

\def\@oval(#1,#2)[#3]{\begingroup\boxmaxdepth \maxdimen
  \@ovttrue \@ovbtrue \@ovltrue \@ovrtrue
  \@tfor\@tempa :=#3\do{\csname @ov\@tempa false\endcsname}\@ovxx
  #1\unitlength \@ovyy #2\unitlength
  \@tempdimb \ifdim \@ovyy >\@ovxx \@ovxx\else \@ovyy \fi
  \advance \@tempdimb -2pt\relax  
  \@getcirc \@tempdimb
  \@ovro \ht\@tempboxa \@ovri \dp\@tempboxa
  \@ovdx\@ovxx \advance\@ovdx -\@tempdima \divide\@ovdx \tw@
  \@ovdy\@ovyy \advance\@ovdy -\@tempdima \divide\@ovdy \tw@
  \@circlefnt \setbox\@tempboxa
  \hbox{\if@ovr \@ovvert32\kern -\@tempdima \fi
  \if@ovl \kern \@ovxx \@ovvert01\kern -\@tempdima \kern -\@ovxx \fi
  \if@ovt \@ovhorz \kern -\@ovxx \fi
  \if@ovb \raise \@ovyy \@ovhorz \fi}\advance\@ovdx\@ovro
  \advance\@ovdy\@ovro \ht\@tempboxa\z@ \dp\@tempboxa\z@
  \@put{-\@ovdx}{-\@ovdy}{\box\@tempboxa}%
  \endgroup}

\def\@ovvert#1#2{\vbox to \@ovyy{%
    \if@ovb \@tempcntb \@tempcnta \advance \@tempcntb by #1\relax
      \kern -\@ovro \hbox{\char \@tempcntb}\nointerlineskip
    \else \kern \@ovri \kern \@ovdy \fi
    \leaders\vrule width \@wholewidth\vfil \nointerlineskip
    \if@ovt \@tempcntb \@tempcnta \advance \@tempcntb by #2\relax
      \hbox{\char \@tempcntb}%
    \else \kern \@ovdy \kern \@ovro \fi}}

\def\@ovhorz{\hbox to \@ovxx{\kern \@ovro
    \if@ovr \else \kern \@ovdx \fi
    \leaders \hrule height \@wholewidth \hfil
    \if@ovl \else \kern \@ovdx \fi
    \kern \@ovri}}

\def\circle{\@ifstar{\@dot}{\@circle}}
\def\@circle#1{\begingroup \boxmaxdepth \maxdimen \@tempdimb #1\unitlength
   \ifdim \@tempdimb >15.5pt\relax \@getcirc\@tempdimb
      \@ovro\ht\@tempboxa 
     \setbox\@tempboxa\hbox{\@circlefnt
      \advance\@tempcnta\tw@ \char \@tempcnta
      \advance\@tempcnta\m@ne \char \@tempcnta \kern -2\@tempdima
      \advance\@tempcnta\tw@
      \raise \@tempdima \hbox{\char\@tempcnta}\raise \@tempdima
        \box\@tempboxa}\ht\@tempboxa\z@ \dp\@tempboxa\z@
      \@put{-\@ovro}{-\@ovro}{\box\@tempboxa}%
   \else  \@circ\@tempdimb{96}\fi\endgroup}

\def\@dot#1{\@tempdimb #1\unitlength \@circ\@tempdimb{112}}

\def\@circ#1#2{\@tempdima #1\relax \advance\@tempdima .5pt\relax
   \@tempcnta\@tempdima \@tempdima 1pt\relax
   \divide\@tempcnta\@tempdima 
   \ifnum\@tempcnta > 15\relax \@tempcnta 15\relax \fi    
   \ifnum \@tempcnta >\z@ \advance\@tempcnta\m@ne\fi
   \advance\@tempcnta #2\relax
   \@circlefnt \char\@tempcnta}

\font\tenln line10
\font\tencirc lcircle10
\font\tenlnw linew10
\font\tencircw lcirclew10

\thinlines   

\newcount\@xarg
\newcount\@yarg
\newcount\@yyarg
\newcount\@multicnt 
\newdimen\@xdim
\newdimen\@ydim
\newbox\@linechar
\newdimen\@linelen
\newdimen\@clnwd
\newdimen\@clnht
\newdimen\@dashdim
\newbox\@dashbox
\newcount\@dashcnt
\catcode`@=12

\catcode`@=11
\font\@linefnt linew10 at 2.4pt
\catcode`@=12

\font\grand=cmr10 at 14pt

\def\u{\underline }
\def\P{{\cal P}}
\def\K{{\cal K}}

\def\g{\mathop{\dashv}\nolimits}

\def\d{\mathop{\vdash}\nolimits}

\def\m{\mathop{\perp}\nolimits}

\def\l{\mathop{\prec}\nolimits}

\def\r{\mathop{\succ}\nolimits}


\def\N{\noindent}
\def\S{\smallskip \par}
\def\M{\medskip \par}
\def\B{\bigskip \par}
\def\BB{\bigskip \bigskip \par}


\def\cc{\gamma}

\def\DD{\Delta}

\def\SS{\Sigma}

\def\sqr#1#2{{\vcenter{\vbox{\hrule height.#2pt
\hbox{\vrule width .#2pt height#1pt \kern#1pt
\vrule width.#2pt}
\hrule height.#2pt}}}}

\def \Vect{\mathop{\rm Vect}\nolimits}

\def\t{\otimes }

\def\tT#1#2{#1^{\otimes #2}}

\def\row#1#2#3{(#1_{#2},\ldots,#1_{#3})}

\def\arbreA{\kern-0.4ex
\hbox{\unitlength=.25pt
\picture(60,40)(0,0)
\put(30,0){\droite(0,1){10}}
\put(30,10){\droite(-1,1){30}}
\put(30,10){\droite(1,1){30}}
\endpicture}\kern 0.4ex}

\def\arbreB{\kern-0.4ex
\hbox{\unitlength=.25pt
\picture(60,40)(0,0)
\put(30,0){\droite(0,1){10}}
\put(30,10){\droite(-1,1){30}}
\put(30,10){\droite(1,1){30}}
\put(15,25){\droite(1,1){15}}
\endpicture}\kern 0.4ex}

\def\arbreC{\kern-0.4ex
\hbox{\unitlength=.25pt
\picture(60,40)(0,0)
\put(30,0){\droite(0,1){10}}
\put(30,10){\droite(-1,1){30}}
\put(30,10){\droite(1,1){30}}
\put(45,25){\droite(-1,1){15}}
\endpicture}\kern 0.4ex}

\def\arbreBC{\kern-0.4ex
\hbox{\unitlength=.25pt
\picture(60,40)(0,0)
\put(30,0){\droite(0,1){40}}
\put(30,10){\droite(-1,1){30}}
\put(30,10){\droite(1,1){30}}
\endpicture}\kern 0.4ex}

\def\arbreun{\kern-0.4ex
\hbox{\unitlength=.25pt
\picture(60,40)(0,0)
\put(30,0){\droite(0,1){20}}
\put(30,20){\droite(-1,1){30}}
\put(30,20){\droite(1,1){30}}
\put(20,30){\droite(1,1){20}}
\put(10,40){\droite(1,1){10}}
\endpicture}\kern 0.4ex}

\def\arbredeux{\kern-0.4ex
\hbox{\unitlength=.25pt
\picture(60,40)(0,0)
\put(30,0){\droite(0,1){20}}
\put(30,20){\droite(-1,1){30}}
\put(30,20){\droite(1,1){30}}
\put(20,30){\droite(1,1){20}}
\put(30,40){\droite(-1,1){10}}
\endpicture}\kern 0.4ex}

\def\arbretrois{\kern-0.4ex
\hbox{\unitlength=.25pt
\picture(60,40)(0,0)
\put(30,0){\droite(0,1){20}}
\put(30,20){\droite(-1,1){30}}
\put(30,20){\droite(1,1){30}}
\put(50,40){\droite(-1,1){10}}
\put(10,40){\droite(1,1){10}}
\endpicture}\kern 0.4ex}

\def\arbrequatre{\kern-0.4ex
\hbox{\unitlength=.25pt
\picture(60,40)(0,0)
\put(30,0){\droite(0,1){20}}
\put(30,20){\droite(-1,1){30}}
\put(30,20){\droite(1,1){30}}
\put(40,30){\droite(-1,1){20}}
\put(30,40){\droite(1,1){10}}
\endpicture}\kern 0.4ex}

\def\arbrecinq{\kern-0.4ex
\hbox{\unitlength=.25pt
\picture(60,40)(0,0)
\put(30,0){\droite(0,1){20}}
\put(30,20){\droite(-1,1){30}}
\put(30,20){\droite(1,1){30}}
\put(40,30){\droite(-1,1){20}}
\put(50,40){\droite(-1,1){10}}
\endpicture}\kern 0.4ex}

\def\simplexes{\kern-0.4ex
\hbox{\unitlength=.25pt
\picture(800,250)(0,0)
\put(0,0){$\Delta^0$}
\put(10,50){$\bullet$}
\put(10,100){$\underline 1$}

\put(200,0){$\Delta^1$}
\put(160,50){$\bullet$}
\put(170,60){\droite(1,0){100}}
\put(260,50){$\bullet$}
\put(160,100){$\underline 1$}
\put(260,100){$\underline 2$}

\put(450,0){$\Delta^2$}
\put(410,50){$\bullet$}
\put(420,60){\droite(1,0){100}}
\put(420,60){\droite(1,2){50}}
\put(520,60){\droite(-1,2){50}}
\put(410,50){$\bullet$}
\put(390,80){$\underline 1$}
\put(460,150){$\bullet$}
\put(530,80){$\underline 2$}
\put(510,50){$\bullet$}
\put(460,180){$\underline 3$}

\endpicture}\kern 0.4ex}

\def\stasheff{\kern-0.4ex
\hbox{\unitlength=.25pt
\picture(800,270)(0,0)
\put(0,0){$\K^0$}
\put(10,50){$\bullet$}
\put(0,100){$\arbreA$}

\put(200,0){$\K^1$}
\put(160,50){$\bullet$}
\put(170,60){\droite(1,0){120}}
\put(280,50){$\bullet$}
\put(140,100){$\arbreB$}
\put(270,100){$\arbreC$}
\put(210,70){$\arbreBC$}

\put(590,0){$\K^2$}
\put(520,130){\droite(1,1){50}}
\put(570,80){\droite(3,1){50}}
\put(620,100){\droite(0,1){60}}
\put(570,180){\droite(3,-1){50}}
\put(520,130){\droite(1,-1){50}}

\put(610,150){$\bullet$}
\put(640,150){$\arbredeux$}

\put(560,170){$\bullet$}
\put(540,200){$\arbreun$}

\put(510,120){$\bullet$}
\put(460,110){$\arbretrois$}

\put(560,70){$\bullet$}
\put(640,70){$\arbrequatre$}

\put(610,90){$\bullet$}
\put(530,20){$\arbrecinq$}

\endpicture}\kern 0.4ex}

\centerline {\grand Une dualit\'e entre simplexes standards et polytopes de Stasheff}
\BB

\centerline {\grand A duality between standard simplices and Stasheff polytopes}
\BB
\centerline {\bf Jean-Louis Loday and Mar\' \i a Ronco}
\BB

\N {\bf R\'esum\'e.}  On montre que la famille des modules de chaines des simplexes standards
peut \^etre munie d'un structure d'op\'erade. De m\^eme la famille des modules de cochaines
des polytopes de Stasheff peut \^etre munie d'un structure d'op\'erade. On montre que, d'une part,
ces deux op\'erades sont duales l'une de l'autre, et, d'autre part, que ce sont des op\'erades de
Koszul.

Les alg\`ebres sur l'op\'erade des simplexes standards, appell\'ees trig\`ebres associatives, sont
d\'etermin\'ees par 3 op\'erations et 11 relations. Les alg\`ebres sur l'op\'erade des polytopes de
Stasheff, appell\'ees trig\`ebres dendriformes, sont d\'etermin\'ees par 3 op\'erations et 7
relations. 
\B

\N {\bf Abstract.}  We show that the family  of chain modules over the standard simplices can be
equipped with an operad structure. Similarly,  the family of  cochain modules of the Stasheff
polytopes can be equipped with an operad structure. We first  show that these operads are Koszul
dual to each other, and second that they are both Koszul operads. 

The algebras over the standard simplices operad, called  {\it associative trialgebras}, are defined
by 3 operations and 11 relations. The algebras over the Stasheff polytopes  operad, called  {\it
dendriform trialgebras}, are defined by 3 operations and 7 relations. 
\B

\N {\bf English version.} Let $\DD^n$ be the the standard simplex of dimension $n$. We show that
the   chain module $\P _{\DD }(n) := C_{*}(\DD^{n-1}), n\geq 1$, can be equipped with a
structure of non-$\SS$-operad (cf. [O]).  Let $\K^n$ be the the
Stasheff polytope of dimension $n$. We show that the   cochain module $\P ^{\K} (n) :=
C^{*}(\K^{n-1}), n\geq 1$, can be equipped with a structure of non-$\SS$-operad. Both operads
are quadratic and we show that they are dual to each other in the operadic sense (cf. [GK]).
Our main result is to show that both operads are Koszul operads, that is : the associated Koszul complexes are
acyclic.

Both operads are binary, generated by three operations, one for each cell of the interval $\DD^1 =
\K^1$. Hence the algebras over $\P_{\DD}$ are determined by two operations of degree 0 and one of
degree 1, and by 11 relations (one for each of the cells of the pentagon $\K^2$). They are called {\it
associative trialgebras}, since all the relations are of the associativity type, cf. 2.  The
algebras over $\P ^{\K} $ are determined by two operations of degree 0 and one of degree 1, and by
7 relations (one for each of the cells of the triangle $\DD^2$). They are called {\it
dendriform trialgebras}, because the cells of the Stasheff polytope can be parametrized by the
planar trees, cf. 3.  
\B

\N {\bf Version fran\c caise.} 
\M
\N {\bf Convention.} On travaille sur un corps $K$ et le produit tensoriel au-dessus de $K$ est
not\' e $\t$. Le produit tensoriel de $n$ copies de l'espace vectoriel V est not\' e $\tT Vn$.
\M

 \N {\bf 1. L'op\'erade des simplexes standards.} Le simplexe standard $\DD^n$ est
par d\'efinition le polytope 
$$\displaylines{
\DD^n := \{ \row x0n  \in {\bf R}^{n+1} \mid  0\leq x_i\leq 1 \hbox { and } x_0 +\cdots +x_n =1 \} :\cr
\simplexes\cr
}$$
Le sommet $\u i$ a toutes ses coordonn\'ees nulles sauf la $i-1$ -i\`eme qui vaut 1. Une face de
dimension $d$ de $\DD^{n-1}$ est compl\`etement d\'etermin\'ee par ses sommets, c'est \`a dire
par un sous-ensemble \`a $d$ \'el\'ements de $[n] := \{ \u 1 , \cdots , \u n \}$. Le complexe de
chaines de $\DD^{n-1}$ est not\'e $C_*(\DD^{n-1})=\oplus_{d\geq 0} C_d(\DD^{n-1})$. On prendra pour base de
$C_d(\DD^{n-1})$ les sous-ensembles \`a $d$ \'el\'ements de $[n] $, car chacun de ces sous-ensembles d\'etermine
une cellule de dimension $d$ de $\DD^{n-1}$.

On d\'efinit une structure d'op\'erade sur les $\P _{\DD }(n) :=
C_*(\DD^{n-1})$ de la mani\`ere suivante. L'application de composition $$
\cc : \P _{\DD }(n)\t \P _{\DD }(i_1)\t \cdots \t \P _{\DD }(i_n) \longrightarrow \P _{\DD }(i_1+
\cdots +i_n)$$
est compl\`etement d\'etermin\'ee par sa valeur sur les vecteurs de base des $\P _{\DD }(k)$.
Notons 
$$
{\rm bij} : [{i_1}]\cup \cdots \cup  [ {i_n}] \to  [{i_1+ \cdots +i_n}]
$$
la bijection \'evidente. Soit $X= \{ \u {j_1},\cdots ,  \u {j_k}\} \subset  [ {n}]; X_1 \subset 
[ {i_1}], \cdots  , X_n \subset [{i_n}]$ des vecteurs de base. On pose
$$
\cc (X; X_1, \cdots , X_n ) := {\rm bij }\,  ( X_{j_1} \cup \cdots \cup X_{j_k}) \subset  [{i_1+ \cdots +i_n}].
$$
\M
\N {\bf 1.1. Proposition.} {\it La composition $\cc$ d\'efinie ci-dessus munit les
 $\P _{\DD }(n), n\geq 1$ d'une structure d'op\'erade (non sym\'etrique).}
\M
Rappelons que l'on travaille dans le cadre des op\'erades non sym\'etriques, donc le foncteur 
 ${\P _{\DD}} : \Vect \to \Vect $ est donn\'e par  $\P _{\DD }(V) = \oplus_{n\geq
1}   \P _{\DD }(n)\t \tT Vn$.
\M

\N {\bf 2. Les trig\`ebres associatives.} Par d\'efinition une {\it trig\`ebre associative} est la
donn\'ee d'un espace vectoriel $A$ et de trois op\'erations  :
$$\displaylines{
\g : A\t A \to A\qquad \hbox {(op\'eration gauche)},\cr
\d : A\t A \to A\qquad  \hbox {(op\'eration droite)},\cr
\m : A\t A \to A\qquad  \hbox {(op\'eration milieu)},\cr
}$$
 satisfaisant aux 11 relations suivantes :
$$\eqalign{
(x\g y)\g z &= x\g (y\g z) , \cr
(x\g y)\g z &= x\g (y\d z) , \cr 
(x\d y)\g z &= x\d (y\g z) , \cr
(x\g y)\d z &= x\d (y\d z) , \cr
(x\d y)\d z &= x\d (y\d z) , \cr
 & \cr
(x\g y)\g z &= x\g (y\m z) , \cr
(x\m y)\g z &= x\m (y\g z) , \cr 
(x\g y)\m z &= x\m (y\d z) , \cr
(x\d y)\m z &= x\d (y\m z) , \cr
(x\m y)\d z &= x\d (y\d z) , \cr
 & \cr
(x\m y) \m z &= x\m (y\m z) . \cr
}$$
\M

\N {\bf 2.1. Th\' eor\`eme.} {\it L'op\' erade (non sym\' etrique) $\P _{\DD }$ des simplexes standards
est binaire et quadratique. Une alg\`ebre sur $\P _{\DD }$ est une trig\`ebre associative. En particulier, la
trig\`ebre associative libre sur un g\'en\'erateur  $\P _{\DD }(K)$ s'identifie \`a 
$\oplus_{n\geq 1} C_*(\Delta ^{n-1})$ muni des op\'erations induites par
$$
X \g Y = {\rm bij }\ (X), \quad X \d Y = {\rm bij }\ (Y), \quad X \m Y = {\rm bij }\ (X\cup Y). $$}
Si  l'on oublie  l'op\' eration $\m$ , alors  $A$ est tout simplement une dig\`ebre associative au sens de
[L1, L2].

Les  alg\`ebres associatives sont des exemples particuliers de trig\`ebre associative en prenant 
$\g = \d = \m$.
\M
Comme l'op\' erade $\P _{\DD }$ est quadratique, elle admet une op\' erade duale au sens de la
dualit\' e de Koszul des op\' erades d\^ue \`a Ginzburg et Kapranov [GK]. Son calcul, et en
particulier son lien avec les polytopes de Stasheff, est l'objet des paragraphes suivants.
\M
\N {\bf 2.2. Trig\`ebres associatives diff\'erentielles.} Une trig\`ebre associative $A$ est dite diff\'erentielle
gradu\'ee, si l'on a $A = \oplus_{n\geq 1} A_n$ et une application lin\'eaire $d : A \to A$ de degr\'e $-1$ et de
carr\'e nul, v\'erifiant de surcroit:
$$\eqalign{
d(x\g y) &= dx \g y , \cr
d(x\d y) &= (-1)^{\mid x \mid} x \d dy , \cr
d(x\m y) &= dx \m y +(-1)^{\mid x \mid} x \m dy . \cr
}$$
La somme de  complexes de chaines $\oplus_{n\geq 1}C_*(\DD^{n-1})$ forme un trig\`ebre associative
diff\'erentielle gradu\'ee.
\B

\N {\bf 3. Les trig\`ebres dendriformes.} Par d\' efinition 
 une {\it trig\`ebre dendriforme} est la donn\'ee d'un espace vectoriel 
$D$
 et de trois op\'erations :
$$\displaylines{
\l : D\t D \to D\qquad \hbox {  (op\'eration gauche)},\cr
\r : D\t D \to D\qquad \hbox {   (op\'eration droite)},\cr
\cdot\,  : D\t D \to D\qquad \hbox {   (op\'eration milieu)},\cr
}$$
satisfaisant aux 7 relations suivantes :
$$\eqalign {
(x \l  y) \l  z &= x \l  (y * z)\ , \cr
(x \r  y) \l  z &= x \r  (y \l  z)\ , \cr
(x * y) \r  z &= x \r  (y \r  z)\ , \cr
\cr
(x \r  y) \cdot  z &= x \r  (y \cdot  z)\ , \cr
(x \l  y) \cdot  z &= x \cdot  (y \r  z)\ , \cr
(x \cdot  y) \l  z &= x \cdot  (y \l  z)\ , \cr
\cr
(x \cdot  y) \cdot  z &= x \cdot  (y \cdot  z)\ . \cr
}$$
o\`u $x*y := x\l y + x\r y + x\cdot y$.
\M

\N {\bf 3.1. Th\'eor\`eme.} {\it L'op\' erade des trig\`ebres dendriformes
 est duale, au sens de la dualit\' e de Koszul des op\' erades, de l'op\'erade 
des trig\`ebres associatives.}
\M
L'op\' erade $\P _{\DD }$ est engendr\' ee par les trois op\' erations $\g, \d, \m$, donc
l'op\'erade duale ${\P _{\DD }}^!$ est aussi engendr\' ee par 3 op\' erations que l'on note $\l, \r,
\cdot$. L'espace des op\' erations que l'on peut faire avec 3 variables est $2\times 3^2 = 18$.
L'espace des relations de ${\P _{\DD }}^!$ est l'orthogonal, pour une certaine forme  quadratique 
non  d\'eg\' en\' er\' ee, de l'espace des relations de $\P _{\DD }$. Il est donc de dimension
18-11=7. On v\' erifie qu'il est bien engendr\' e par les relations de d\' efinition des trig\`ebres
dendriformes.
\M

Si l'op\'eration $\cdot$ est triviale, alors la trig\`ebre dendriforme  $D$ est  une dig\`ebre
dendriforme au sens de [L1, L2].

On remarque ais\'ement, en ajoutant toutes les relations, que le produit $*$, somme des trois
produits $\l, \r, \cdot$, munit $D$ d'une structure d'alg\`ebre associative. Ainsi toute trig\`ebre dendriforme
d\'etermine une alg\`ebre associative. Ce foncteur est dual du foncteur qui va des alg\`ebres associatives dans les
trig\`ebres associatives (rappellons que l'op\'erade des alg\`ebres associatives est auto-duale).
\M

\N {\bf 3.2. L'op\' erade des polytopes de Stasheff.} Le polytope de Stasheff $\K^n$, appel\' e parfois
{\it associa\`edre}, est un polytope simple de dimension $n$, hom\' eomorphe \`a la boule unit\' e
et dont les sommets sont en bijection avec les parenth\' esages associatifs d'un mot \`a $n+2$
lettres, cf. [St]. De mani\`ere \' equivalente on peut identifier ces mots parenth\'es\'es aux (classes
d'isotopie d') arbres binaires planaires \`a $n+2$ feuilles (et donc $n+1$ sommets internes). Les
faces de dimension $d$ de $\K^n$ sont en bijection avec les arbres planaires \`a  $n+2$ feuilles
et  $n+1-d$ sommets internes. On note $C^*(\K^n)= \oplus_{d\geq 0} C^d(\K^n)$ le module de cochaines de $\K^n$
et on prend pour base de $C^d(\K^n)$ les arbres planaires \`a  $n+2$ feuilles
et  $n+1-d$ sommets internes :
$$\stasheff$$
Rappelons que tout arbre planaire \`a $n+1$ feuilles $x$ est le  greff\' e d'un certain nombre d'arbres : 
$x=x^{(1)}\vee \cdots \vee x^{(k)}$ avec $k\geq 2$ si
$x\not= \vert$.
\M

\N {\bf 3.3. Th\' eor\`eme.} {\it L'op\' erade binaire quadratique (non sym\' etrique) $\P^{\K}$
des trig\`ebres dendriformes est telle que 
$$
\P ^{\K}(n) = C^*(\K^{n-1}).
$$
En particulier, la
trig\`ebre dendriforme libre sur un g\'en\'erateur  $\P^{\K }(K)$ s'identifie \`a 
$\oplus_{n\geq 1} C^*(\K ^{n-1})$ muni des op\'erations induites par
$$\eqalign{
x\l y &= x^{(1)}\vee \cdots \vee ( x^{(k)}* y)\, , \cr
x\r y &= (x*y^{(1)})\vee \cdots \vee y^{(\ell)} \, , \cr
x\cdot y &= x^{(1)}\vee \cdots \vee (x^{(k)}* y^{(1)})\vee \cdots \vee y^{(\ell)}\,  ,\cr
}$$}
o\`u $y=y^{(1)}\vee \cdots \vee y^{(\ell)}$.

Dans ces formules r\' ecursives l'arbre $\vert$ sans sommet interne est \'el\'ement neutre pour l'op\' eration $*$\ .
\M

 \N {\bf 3.4. Remarque.} Il est bien connu que les complexes de chaines $C_*(\K^{n})$ forment une op\' erade
(r\'egissant les
$A_{\infty}$-alg\`ebres), cf [St]. Mais dans ce cas $C_*(\K^{n})$ est en degr\' e $n+2$ alors que dans notre
contexte il est en degr\'e $n+1$. On a donc \`a faire \`a une tout autre op\' erade.
\M
\N {\bf 3.5. Compatibilit\'e avec la diff\'erentielle.}  Les espaces vectoriels $\P ^{\K}(n) = C^*(\K^{n-1})$ sont
filtr\'es par la dimension des cellules. La structure d'op\'erade est compatible avec cette filtration. L'op\'erade
gradu\'ee associ\'ee est celle construite par Chapoton dans [Ch].
\M
\N {\bf 3.6. Comultiplication.} La dig\`ebre (resp. trig\`ebre) dendriforme libre est naturellement munie d'une
comultiplication dont l'\'etude est faite dans [R] (resp. dans une publication ult\'erieure).
\B

\N {\bf 4. Op\' erades de Koszul.} A toute op\' erade quadratique $\P$ ou peut associer un
complexe de Koszul (cf. [GK]). Lorsque celui-ci est acyclique on dit que l'op\' erade est de
Koszul. L'une des cons\' equences de cette propri\' et\' e est l'existence d'un ``petit" complexe
explicite, construit \`a partir de l'op\' erade duale $\P ^!$  pour calculer l'homologie de
Quillen d'une $\P$-alg\`ebre $A$. On note $C^{\P}_*(A)$ ce complexe et $H^{\P}_*(A)$ son
homologie.
\M
\N {\bf 4.1. Th\' eor\`eme.} {\it Le complexe de chaines d'une trig\`ebre dendriforme $A$ est donn\'e par
$$C^{\P^{\K}}_n(A)= C_*(\Delta^{n-1})\t A^{\t n} , \quad d = - \sum_{i=1}^{i=n-1} (-1)^i d_i ,$$
avec $d_i(X; a_1,\cdots ,a_n) = (d_i(X); a_1,\cdots ,a_i \circ_i^X a_{i+1}, \cdots a_n),$
o\`u $d_i(X)$ est l'image de $X$ par l'application $d_i : [n] \to [n-1]$ donn\'ee par 
$$
d_i({\u r}) = \cases{
\u {r-1}  &si  $ i\leq r-1,$\cr
{\u r}  & si $ i\geq r.$\cr} 
$$
et o\`u $ \circ_i^X$ est donn\'e par
$$
\circ_i^X = \cases{
\cdot    &  si $ i\in X \hbox { and } i+1 \in X,$\cr
\prec    & si $ i\notin X \hbox { and } i+1 \in X,$\cr
\succ    & si $ i\in X \hbox { and } i+1 \notin X,$\cr
*          & si $ i\notin X \hbox { and } i+1 \notin X.$\cr}
$$}
\N {\bf 4.2. Th\' eor\`eme.} {\it Le complexe de chaines d'une trig\`ebre associative $A$ est donn\'e par
$$C^{\P _{\DD}}_n(A)= C_*(\K^{n-1})\t A^{\t n} , \quad d = - \sum_{i=1}^{i=n-1} (-1)^i d_i ,$$
avec $d_i(t; a_1,\cdots ,a_n) = (d_i(t); a_1,\cdots ,a_i \circ_i^t a_{i+1}, \cdots a_n),$
o\`u $d_i(t)$ est l'arbre obtenu \`a partir de $t$ en supprimant la $i$-i\`eme feuille, 
et o\`u $ \circ_i^t$ est donn\'e par
$$
\circ_i^t = \cases{
\g   & si $ \hbox { la $i$-\`eme feuille de $t$ est orient\'ee vers la gauche} ,$\cr
\d    & si  $ \hbox { la $i$-\`eme feuille de $t$ est orient\'ee vers la droite} ,$\cr
\m    & si $ \hbox { la $i$-\`eme feuille de $t$ est l'une des feuilles du milieu} .$\cr}
$$}

\N {\bf 4.3. Th\' eor\`eme.} { \it Les deux op\' erades $\P _{\DD}$ et $\P ^{\K}$ sont des op\' erades de
Koszul.}
\M
Puisqu'elles sont duales l'une de l'autre il suffit de montrer que $\P _{\DD}$ est de Koszul. On
sait, d'apr\`es [GK],  qu'il est \' equivalent de montrer que l'homologie de la $\P$-alg\`ebre libre $\P (V)$ est
triviale, plus pr\' ecis\' ement que 
$$ H^{\P}_1(\P (V) )= V \quad {\rm et } \quad H^{\P}_n(\P (V))= 0 \ {\rm pour }\ n>1.$$
En fait il suffit de regarder le cas $V=K$. On montre alors que le complexe 
 $C^{\P^{\K}}_*( \P  _{\DD}(K) )$ est somme directe (index\'ee par les cellules des simplexes standards) de
complexes de chaines augment\'es de certains ensembles simpliciaux. On ach\`eve la preuve en
montrant que ces ensembles simpliciaux sont  contractiles.
\M

\N {\bf 4.4. S\'eries g\'en\'eratrices.} Lorsqu'une op\' erade quadratique non sym\' etrique $\P$ est de Koszul,
d'op\'erade duale $\P^!$,  sa s\'
erie g\' en\' eratrice 
$$f^{\P }(x) := \sum_{n\geq 1} (-1)^n \dim \P(n) x^n$$
 v\' erifie la propri\' et\' e suivante, cf.
[GK]:  $f^{\P}(f^{\P^!}(x))=x$.
On peut affiner ce r\' esultat  dans le cas des op\' erades quadratiques \`a
valeurs dans les espaces vectoriels filtr\' es en rempla\c cant $\dim \P(n)$ par le polynome de
Poincar\' e $p(n,t) = \sum_{i\geq 0}p_i(n) t^i$, o\`u $p_i(n)$ est la dimension de l'espace
$F^i\P(n)/F^{i-1}\P(n) $. On obtient alors la s\' erie g\' en\' eratrice  
$f^{P }_t(x) := \sum_{n\geq 1} (-1)^n p(n,t) x^n$. Lorsque $\P$ est de Koszul les s\' eries $ 
f^{\P}_t$ et $f^{\P^!}_t$ sont encore inverses l'une de l'autre.

Pour $\P _{\DD}$ on obtient $p(n,t) = {1 \over t} ((1+t)^n - 1)$ et donc 
$$f^{\DD}_t(x) = {-x \over (1+x)(1+(1+t)x)}\ .$$
En cons\'equence du th\'eor\`eme 4.1,  $f^{\K}_t(x)$ est la s\'erie g\'en\'eratrice inverse.  On obtient
$$f^{\K}_t(x) =  {-(1+(2+t)x) +\sqrt { 1+2(2+t)x+t^2x^2
}\over 2(1+t)x}\ . $$
On a ainsi obtenu, comme sous-produit, la s\'erie g\'en\'eratrice des arbres planaires. On v\'erifie
imm\'ediatement que  $f^{\K}_0$ est la s\'erie g\'en\'eratrice des nombres de Catalan et  $f^{\K}_1$  la s\'erie
g\'en\'eratrice des super nombres de Catalan.
\B

\N {\bf 5. Op\' erades cubiques.} Une dualit\' e similaire existe pour la famille de polytopes
$I^n$, c'est \`a dire les hypercubes. Les complexes de modules de chaines $\P_Q(n):= C_*(I^{n-1})$
s'assemblent pour donner une op\' erade. L'op\' erade
duale est donn\' ee par  $\P^Q(n):= C^*(I^{n-1})$. Les alg\`ebres associ\' ees sont les
{\it trig\`ebres cubiques} d\' efinies par 3 op\' erations et 9 relations. Ces deux op\' erades sont
de Koszul et leur s\' erie g\' en\' eratrice commune est 
$f^{I}_t(x) = {-x \over 1+(t+2)x}\ .$
On v\' erifie imm\' ediatement que l'on a:
$ f^{I}_t(f^{I}_t(x) ) = x\ .$
\B

\centerline{\bf Ref\' erences}
\B
\N [Ch] Chapoton Fr\'ed\'eric, {\it Op\'erades, polytopes et big\`ebres}, th\`ese, Universit\'e Paris
VI, 2000.
\S
\N [GK]  Ginzburg, Victor; Kapranov, Mikhail. {\it  Koszul duality for operads}. Duke Math. J. 76
(1994), no. 1, 203--272.
\S
\N [L1]  Loday, Jean-Louis. {\it  Alg\`ebres ayant deux op\'erations associatives (dig\`ebres)}. C.
R. Acad. Sci. Paris S\'er. I Math. 321 (1995), no. 2, 141--146. 
\S
\N [L2]   Loday, Jean-Louis. {\it  Dialgebras}, preprint IRMA Strasbourg, 1999.
\S
\N [LR1] Loday, Jean-Louis; Ronco, Mar\' \i a O. {\it Hopf algebra of the planar binary trees}. Adv.
Math. 139 (1998), no. 2, 293--309.
\S
\N [LR2]  Loday, Jean-Louis; Ronco, Mar\' \i a O. {\it Trialgebras}, en pr\'eparation.
\S
\N [O] Operads: Proceedings of
Renaissance Conferences (Hartford, CT/Luminy, 1995), 37--52, Contemp. Math., 202, Amer. Math. Soc., Providence, RI, 1997.
\S
\N [R] Ronco, Mar\' \i a O. {\it A Milnor-Moore theorem for dendriform Hopf algebras}.  C.
R. Acad. Sci. Paris S\'er. I Math.  331 (2000),
\S
\N [St] Stasheff, James Dillon. {\it  Homotopy associativity of $H$-spaces}. I, II. Trans. Amer.
Math. Soc. 108 (1963), 275-292; ibid. 108 1963 293--312.
\B
\N JLL : Institut de Recherche Math\'ematique Avanc\'ee,

    CNRS et Universit\'e Louis Pasteur, 7 rue R. Descartes,

    67084 Strasbourg Cedex, France

    E-mail : loday@math.u-strasbg.fr
\M

\N MOR : Departamento de Matem\'atica

Ciclo B\'asico Com\'un, Universidad de Buenos Aires

Pab. 3 Ciudad Universitaria Nu\~nez 

(1428) Buenos-Aires, Argentina

E-mail : mronco@mate.dm.uba.ar
\B

\hfill PolytopesDualiteNote 29.01.01

\end